\documentclass{amsart}
\usepackage{amsmath,amsthm,amssymb}

\usepackage{lineno,hyperref}
\usepackage{caption}
\usepackage{subcaption}
\usepackage{setspace}
\usepackage{amsmath, amsthm, amssymb, graphicx}
 \newtheorem{thm}{Theorem}[section]

 \newtheorem{test}[thm]{Test}

 \newtheorem{defn}[thm]{Definition}

 \numberwithin{equation}{section}

\numberwithin{equation}{section}
\begin{document}

\title{Fuzzy general linear methods}

\author{Javad Farzi,
        Afsaneh Moradi, \\ Department of Mathematics, Sahand University Of Technology, P.O. Box 51335-1996,  Tabriz, IRAN.}
        \email[Corresponding author]{farzi@sut.ac.ir}



\maketitle

\begin{abstract}

This paper concerns with the developing the most general schemes
so-called Fuzzy General Linear Methods (FGLM) for solving fuzzy
differential equations. The general linear methods (GLM) for
ordinary differential equations are the middle state of two extreme
extensions (linear multistep and Runge-Kutta methods) of the one
step Euler method. In this paper we develop the FGLM framework of
the Adams schemes for solving fuzzy differential equations under the
strongly generalized differentiability. The stability, consistency
and convergent results will be addressed. The numerical results and
the order of accuracy is illustrated to show the efficiency and
accuracy of the novel scheme.
\end{abstract}

\keywords{   General linear methods, Adams methods, Strongly generalized
differentiability, Generalized fuzzy derivative, Fuzzy differential
equations.}

\subjclass{34A07}

\section{Introduction}
Many problems in science and engineering have some uncertainty in
their nature and fuzzy differential equations are appropriate tools
for modeling of such problems \cite{F1}. The interpretation of a
fuzzy differential equation in the sense of generalized
differentiability allows to fuzzifying the appropriate numerical
methods of ordinary differential equations to fuzzy differential
equations. The Hukuhara derivative with the extension principle or
differential inclusions have some disadvantages. The main drawback
is that the solutions obtained in this setting have increasing
length of their supports \cite{38,Bede3}. Many authors have been
generalized the traditional methods such as Euler's method,
Adams-Bashforth methods, predictor-corrector method, Runge-Kutta
method,...\cite{Bede5,21,20,26,28,29,23,19} to fuzzy differential
and fuzzy initial value problems. However, they use Hukuhara
differentiability and fuzify the numerical method using extension
principle or other methods \cite{F1}. Under the concept of strongly
generalized differentiability there exist fuzzy derivative for a
large class of fuzzy-number-valued functions \cite{Bede2,Bede3}.
Another advantage is that there exist two local solutions,
so-called, (i)-differentiable and (ii)-differentiable solutions.
According to the nature of the initial value problem we can choose
the best meaningful practical solution. In this paper we develop the
GLM schemes based on strongly generalized differentiability concept.
Notion of a fuzzy derivative first introduced by Chang and Zadeh
\cite{3} and Dubo and Prade \cite{4} introduced its extension.
Stefanini \cite{30, 31} introduced the fuzzy gH-difference and Bede
and Stefanini\cite{Bede7} defined and studied new generalization of
the differentiability for fuzzy-number-valued functions. The aim of
this paper is to develop the GLMs for fuzzy differential equations
and study their consistency, stability and convergence. In this
paper, under the strongly generalized differentiability we develop a
well-known Adams-Bashforth methods in the framework of a general
linear method. This starting step will motivate us to develop the
arbitrary classes of GLMs with demanded properties in forthcoming
research.

Let us denote by $\mathbb{R}_{\mathcal{F}}$ the class of fuzzy
numbers, i.e. normal, convex, upper semicontinuous and compactly
supported fuzzy subsets of the real numbers. The fuzzy initial value
problem is defined as follow:
\begin{eqnarray}\label{eq-f}
y'(t) &=& f(t,y(t)),\quad t\in [t_0,T],\\
y(t_0) &=& y_0,\nonumber
\end{eqnarray}
where, $f:[t_0,T]\times \mathbb{R}_{\mathcal{F}}\to
\mathbb{R}_{\mathcal{F}}$ and $y_0\in \mathbb{R}_{\mathcal{F}}$.

Here, we explain the GLM for ordinary IVP (\ref{eq-f}) and in next
sections we will discuss on the development of GLM for FIVP. Burage
and Butcher \cite{2} have presented a standard representation of a
GLM in terms of four matrices. These methods were formulated as
follows:
\begin{equation}\label{1}
\begin{array}{c} Y=hAf(Y)+Uy^{[n-1]},\\
y^{[n]}=hBf(Y)+Vy^{[n-1]}.
\end{array}
\end{equation}
where $y^{[n-1]}$ and $y^{[n]}$ are input and output approximations,
respectively, and
\[A\in \mathbb{R}^{s\times s},\quad U\in \mathbb{R}^{s\times r},\quad B\in \mathbb{R}^{r\times s},\quad V\in \mathbb{R}^{r\times r}.\]
In this paper we use the fuzzy interpolation for constructing
Adams-Bashforth schemes in the general linear methods framework. The
organization of this paper is as follow: In section \ref{sec2} we
present the preliminaries from GLM and fuzzy calculus. In section
\ref{sec3} we apply the GLM form of linear multistep methods to
solve the fuzzy differential equations and in section \ref{sec5}
numerical results are given.

\section{Preliminaries}\label{sec2}
In this section we present the required concepts from general linear
methods and also we shortly review the required definitions form
fuzzy calculus, as given in \cite{Bede1}. We will give the main idea
of the paper for an important subclass of LMMs, the so-called Adams
methods, in GLM framework.

\begin{defn}
Let $u,v\in \mathbb{R}_\mathcal{F}$, the Hukuhara difference
(H-difference $\circleddash_{H}$) of $u$ and $v$ is defined by
\[u\circleddash v=w \Longleftrightarrow u=v+w.\]
\end{defn}
Where $w\in\mathbb{R}_\mathcal{F}$ is called the H-difference of $u$
and $v$. If H-difference $u\circleddash v$ exists, then
$[u\circleddash v]_r=[u_r^--v_r^-,u_r^+-v_r^+]$. The Hukuhara
derivative for a fuzzy function was introduced by Puri and Relescu
\cite{36}.
 From Kaleva \cite{37} and Diamond \cite{38},
it follows that a Hukuhara differentiable function has increasing
length of its support interval. So the Hukuhara difference rarely
exists and to overcome this situation strongly generalized
differentiability of fuzzy-number-valued functions was introduced
and studied by Bede-Gal \cite{Bede3}. Thus, in this case a
differentiable function may have the property that the support has
increasing or decreasing length.
\begin{defn}\label{def2.5}
Let $f:(a,b)\rightarrow\mathbb{R}_\mathcal{F}$ and $x_0\in(a,b)$. We
say that f is strongly generalized differentiable at $x_0$, if there
exists an element $f'(x_0)\in \mathbb{R}_\mathcal{F}$, such that
\begin{itemize}
\item[(i)] for each $h>0$ sufficiently close to 0, the H-differences
$f(x_0+h)\circleddash f(x_0)$ and $f(x_0)\circleddash f(x_0-h)$
exist and
\[\lim_{h\rightarrow0}\frac{f(x_0+h)\circleddash
f(x_0)}{h}=\lim_{h\rightarrow0}\frac{f(x_0)\circleddash
f(x_{0}-h)}{h}=f'(x_0),\] or
\item[(ii)] for each $h>0$ sufficiently close to 0, the H-differences
$f(x_0)\circleddash f(x_0+h) $ and $f(x_0-h)\circleddash f(x_0)$
exist and
\[\lim_{h\rightarrow0}\frac{f(x_0)\circleddash f(x_0+h)}{(-h)}=\lim_{h\rightarrow0}\frac{f(x_0-h)
\circleddash f(x_0)}{(-h)}=f'(x_0).\]
\end{itemize}
\end{defn}
 Let $f:(a,b)\rightarrow\mathbb{R}_\mathcal{F}$, we say that $f$
is $(i)$-differentiable and ($ii$)-differentiable on $(a,b)$ if $f$
is differentiable in the sense ($i$) and ($ii$) of Definition
\ref{def2.5}, respectively. There is also two other
differentiability cases - (iii) and (iv) - differentiability - that
in these cases there is no existence theorems and we do not discuss
them here.

Bede in \cite{Bede5} proved that under certain conditions the fuzzy
initial value problem \eqref{eq-f} has a unique solution and is
equivalent to the system of ODEs
\begin{equation*}
\left\{\begin{array}{c} (y_{r}^-)'=f_{r}^-(t,y_{r}^-,y_r^+)\\
(y_{r}^+)'=f_{r}^+(t,y_{r}^-,y_r^+)
\end{array},r\in[0,1]\right.
\end{equation*}
with respect to H-differentiability.

In this interpretation solutions of a fuzzy differential equation
have always an increasing length of its support interval. So a fuzzy
dynamical system will have more uncertain behavior in time and it
does not allow to have a periodic solutions. Thus, for solve FDEs
the different ideas and methods have been investigated. The second
interpretation was based on Zadeh's extension principle defined in
\cite{44}. Consider the classical ODE $x'=f(t,x,a)$,
$x(t_0)=x_0\in\mathbb{R}$ where $a\in\mathbb{R}$ is a parameter. By
using Zadeh's extension principle on the classical solution, we
obtain a solution of the FIVP. The third interpretation have been
developed based on generalized fuzzy derivative. In this work we
will work with interpretation based on strongly generalized
differentiability. Fuzzy differential equations based on generalized
H-differentiability were investigated by Bede-Gal in \cite{Bede3}
and more general results were proposed in Bede-Gal \cite{Bede4}.
According to the assumptions of the Theorem 9.11 in \cite{Bede1},
the fuzzy initial value problem \eqref{1} is equivalent to the union
of the ODEs:

\begin{eqnarray}\label{i-diff}
\left\{
\begin{array}{ll}
(y_{\alpha}^{-})'(t) =
f_{\alpha}^{-}(t,y_{\alpha}^{-}(t),y_{\alpha}^{+}(t))
\\
(y_{\alpha}^{+})'(t) =
f_{\alpha}^{+}(t,y_{\alpha}^{-}(t),y_{\alpha}^{+}(t)), &\alpha\in
[0,1]
\\
(y_{\alpha}^{-})(t_0) = (y_0)_{\alpha}^{-},\quad
(y_{\alpha}^{+})(t_0) = (y_0)_{\alpha}^{+}.
\end{array}
\right.
\end{eqnarray}
and
\begin{eqnarray}\label{ii-diff}
\left\{
\begin{array}{ll}
(y_{\alpha}^{-})'(t) =
f_{\alpha}^{+}(t,y_{\alpha}^{-}(t),y_{\alpha}^{+}(t))
\\
(y_{\alpha}^{+})'(t) =
f_{\alpha}^{-}(t,y_{\alpha}^{-}(t),y_{\alpha}^{+}(t)), & \alpha\in
[0,1]
\\
(y_{\alpha}^{-})(t_0) = (y_0)_{\alpha}^{-},\quad
(y_{\alpha}^{+})(t_0) = (y_0)_{\alpha}^{+}.
\end{array}
\right.
\end{eqnarray}
For triangular input data we have the same systems \eqref{i-diff}
and \eqref{ii-diff} with an extra equation $(y_{\alpha}^{1})'(t) =
f_{\alpha}^{1}(t,y_{\alpha}^{-}(t),y_{\alpha}^{1}(t),y_{\alpha}^{+}(t))$
where $f=(f^{-},f^{1},f^{+})$ (see Theorem 9.12 in \cite{Bede1}).

A linear multistep method is defined by the first characteristic
polynomial $\rho(r) = \sum_{j=0}^k \alpha_j r^j$ and the second
characteristic polynomial $\sigma(r) = \sum_{j=0}^k \beta_j r^j$ as
follow
\begin{equation}
\sum_{j=0}^k \alpha_j y_{n+j} = h\sum_{j=0}^k \beta_j f_{n+j},
\end{equation}
where $a=t_{n}\leq t_{n+1}\leq \cdots \leq t_{N}=b$,
$h=\frac{b-a}{N}=t_{n+k}-t_{n+k-1}$, $f_{n+j} = f(t_{n+j},y_{n+j})$
and $\alpha_j$ and $\beta_{j}, j=0,1,\cdots,k$ are constants. In
this scheme we can evaluate an approximate solution $y_{n+k}$ for
the exact value $y(x_{n+k})$ using the starting values
$y_0,y_1,\dots, y_{n+k-1}$. The Adams schemes are characterized by
their first characteristic polynomial as $\rho(r)=r^{k}-r^{k-1}$.
Therefor, we have
\begin{equation}\label{eq2.1}
y_{n+k} = y_{n+k-1}+h\sum_{j=0}^{k}\beta_{j}f_{n+j},
\end{equation}
In (\ref{eq2.1}) the case $\beta_{k}=0$ means that the method is
explicit and otherwise the method is implicit. The stability issue
of LMMs are characterized by the root condition for the first
characteristic polynomial $\rho(r)$, that means the roots $r_s,
s=1,2,\dots,k$ of $\rho(r)$ satisfy $|r_s|\le 1$ and the roots with
$|r_s|=1$ are simple \cite{15}. The zero-stability of an LMM and correspondingly the its GLM form depends on that
the first characteristic polynomial $\rho(r)$ or the minimal polynomial of the matrix $V$ satisfies the root condition.

\section{A GLM scheme with strongly generalized differentiability}\label{sec3}
In this section we present the derivation of a GLM based on linear
$k-$step Adams schemes for solving fuzzy initial value problem under
strongly generalized differentiability. Assume that for an equally
spaced points $0=t_0<t_1<\cdot<t_N=T$ at $t_n$ the exact solutions
are indicated by
${\textbf{Y}}_{1}(t_{n};r)=[\textbf{Y}_1^-(t_n;r),\textbf{Y}_1^+(t_n;r)]$
and
$\textbf{Y}_{2}(t_{n};r)=[\textbf{Y}_2^-(t_n;r),\textbf{Y}_2^+(t_n;r)]$
under (i) and (ii)-differentiability, respectively. Also assume that
$y_{1}(t_{n};r)=[y_1^-(t_n;r),y_1^+(t_n;r)]$ and
$y_{2}(t_{n};r)=[y_2^-(t_n;r),y_2^+(t_n;r)]$ are approximate
solutions at $t_n$ under (i) and (ii)-differentiability,
respectively.

The $k$-step Adams methods under Hukuhara or (i)-differentiability
can be written as:
 \begin{equation}\label{equ3.1}
 \begin{array}{ccc}
  y_{1r}^-(t_{n+k};r) &=& y_{1r}^-(t_{n+k-1};r) +h\sum_{j=0}^{k}\beta_j f^-(t_{n+j},y_{1r}(t_{n+j};r)),\\
   y_{1r}^+(t_{n+k};r) &=& y_{1r}^+(t_{n+k-1};r) +h\sum_{j=0}^{k}\beta_j f^+(t_{n+j},y_{1r}(t_{n+j};r)),
 \end{array}
 \end{equation}
 and under (ii)-differentiability can be written as:
  \begin{equation}\label{equ3.2}
  \begin{array}{ccc}
   y_{2r}^-(t_{n+k};r) &=& y_{2r}^-(t_{n+k-1};r) +h\sum_{j=0}^{k}\beta_j f^+(t_{n+j},y_{2r}(t_{n+j};r)),\\
   y_{2r}^+(t_{n+k};r) &=& y_{2r}^+(t_{n+k-1};r) +h\sum_{j=0}^{k}\beta_j f^-(t_{n+j},y_{2r}(t_{n+j};r)),
  \end{array}
 \end{equation}

The Adams schemes are $k$-step methods \eqref{eq2.1} with
$\rho(r)=r^{k}-r^{k-1}$. In this setting we can find their
corresponding general linear method framework. In GLM representation
we should first determine the input and output vectors and then find
the corresponding matrices. For this end we consider the input and
output approximation of general linear methods as follow
\begin{equation*}
y^{[n-1]}=\left(
\begin{array}{c}
y_{n+k-1}\\
hf_{n+k-1}\\
hf_{n+k-2}\\
\vdots\\
hf_{n+1}\\
hf_{n}
\end{array}\right),\qquad y^{[n]}=\left(
\begin{array}{c}
y_{n+k}\\
hf_{n+k}\\
hf_{n+k-1}\\
\vdots\\
hf_{n+2}\\
hf_{n+1}
\end{array}\right).
\end{equation*}
Similarly, a linear k-steps methods under strongly generalized
differentiability \eqref{equ3.1} and \eqref{equ3.2} can be
representation in the form of general linear methods. For this
representation the input vectors for the GLM form of \eqref{equ3.1}
and \eqref{equ3.2} are indicated by
$y_{1r}^{[n-1]}=\big[y_{1r}^{-[n-1]},y_{1r}^{+[n-1]}\big]$ and
$y_{2r}^{[n-1]}=\big[y_{2r}^{-[n-1]},y_{2r}^{+[n-1]}\big]$ under (i)
and (ii)-differentiability, respectively. Corresponding to the input
vectors, the output vectors are indicated by
$y_{1r}^{[n]}=\big[y_{1r}^{-[n]},y_{1r}^{+[n]}\big]$ and
$y_{2r}^{[n]}=\big[y_{2r}^{-[n]},y_{2r}^{+[n]}\big]$ under (i) and
(ii)-differentiability, respectively. Now, we consider the input
approximation of general linear methods in terms of
(i)-differentiability as:
\begin{equation}\label{input_i}
  y_{1r}^{-[n-1]}=\left(\begin{array}{c}
                          y^-_{{n+k-1}_{1r}} \\
                          hf^-_{{n+k-1}_{1r}} \\
                          hf^-_{{n+k-2}_{1r}} \\
                          \vdots\\
                          hf^-_{{n+1}_{1r}}\\
                          hf^-_{{n}_{1r}}
                        \end{array}\right),\qquad
                 y_{1r}^{+[n-1]}=\left(\begin{array}{c}
                          y^+_{{n+k-1}_{1r}} \\
                          hf^+_{{n+k-1}_{1r}} \\
                          hf^+_{{n+k-2}_{1r}} \\
                          \vdots\\
                          hf^+_{{n+1}_{1r}}\\
                          hf^+_{{n}_{1r}}
                        \end{array}\right),
\end{equation}
and under the (ii)-differentiability we obtain the following input
vectors:
 \begin{equation}\label{input_ii}
  y_{2r}^{-[n-1]}=\left(\begin{array}{c}
                          y^-_{{n+k-1}_{2r}} \\
                          hf^+_{{n+k-1}_{2r}} \\
                          hf^+_{{n+k-2}_{2r}} \\
                          \vdots\\
                          hf^+_{{n+1}_{2r}}\\
                          hf^+_{{n}_{2r}}
                        \end{array}\right),\qquad
                 y_{1r}^{+[n-1]}=\left(\begin{array}{c}
                          y^+_{{n+k-1}_{2r}} \\
                          hf^-_{{n+k-1}_{2r}} \\
                          hf^-_{{n+k-2}_{2r}} \\
                          \vdots\\
                          hf^-_{{n+1}_{2r}}\\
                          hf^-_{{n}_{2r}}
                        \end{array}\right),
\end{equation}
By considering the above input vectors, the fuzzy general linear methods form of \eqref{equ3.1} and \eqref{equ3.2} can be formulated in case of (i)-differentiability as:
\begin{equation}\label{GLM1}
  \left(\begin{array}{c}
          Y_{1r} \\ \hline
          y_{1r}^{[n]}
        \end{array}\right)
        =\left(\begin{tabular}{c|c}
                             A & U \\ \hline
                             B & V
                           \end{tabular}\right)\left(\begin{array}{c}
                                                       hf_{1r}(Y_{1r}) \\ \hline
                                                       y_{1r}^{[n-1]}
                                                     \end{array}\right),
\end{equation}
and in case of (ii)-differentiability we have:
\begin{equation}\label{GLM1}
  \left(\begin{array}{c}
          Y_{2r} \\ \hline
          y_{2r}^{[n]}
        \end{array}\right)
        =\left(\begin{tabular}{c|c}
                             A & U \\ \hline
                             B & V
                           \end{tabular}\right)\left(\begin{array}{c}
                                                       hf_{2r}(Y_{2r}) \\ \hline
                                                       y_{2r}^{[n-1]}
                                                     \end{array}\right),
\end{equation}
where $Y_{1r}=[Y_{1r}^-,Y_{1r}^+]$ and $Y_{2r}=[Y_{2r}^-,Y_{2r}^+]$ are internal stages under (i) and (ii)-differentiability, respectively. Also
\begin{equation*}
  \left(\begin{tabular}{c|c}
                             A & U \\ \hline
                             B & V
                           \end{tabular}\right)=
                           \left(\begin{tabular}{c|ccccc}
                                   0 & 1 & $\beta_{k-1}$ & $\cdots$ & $\beta_{1}$ & $\beta_{0}$ \\ \hline
                                   0 & 1 & $\beta_{k-1}$ &$ \cdots$ & $\beta_{1}$ & $\beta_{0}$ \\
                                   1 & 0 & 0 & $\cdots$ & 0 & 0 \\
                                   0 & 0 & 1 & $\cdots$ & 0 & 0 \\
                                   $\vdots$ & $\vdots$ &$ \vdots$ & $\quad$ & $\vdots$ & $\vdots$ \\
                                   0 & 0 & 0 & $\cdots$ & 1 & 0
                                 \end{tabular}\right).
\end{equation*}

Now, we consider two example of Fuzzy GLMs form of k-step methods
under strongly generalized differentiability for $k=4,5$. First,
Consider $k=4$. The input vectors for $k=4$ under (i) and
(ii)-differentiability are as follow, respectively:
\begin{equation*}
{y}^{\mp[n-1]}_{1r}=\left(\begin{array}{c}
{y}^{\mp}_{1r}(t_{n+3})\\
h{f}^{\mp}_{1r}(t_{n+3},y_{1r}(t_{n+3}))\\
h{f}^{\mp}_{1r}(t_{n+2},y_{1r}(t_{n+2}))\\
h{f}^{\mp}_{1r}(t_{n+1},y_{1r}(t_{n+1}))\\
h{f}^{\mp}_{1r}(t_{n},y_{1r}(t_{n}))
\end{array}\right),\quad
{y}^{\mp[n-1]}_{2r}=\left(\begin{array}{c}
{y}^{\mp}_{2r}(t_{n+3})\\
h{f}^{\pm}_{2r}(t_{n+3},y_{2r}(t_{n+3}))\\
h{f}^{\pm}_{2r}(t_{n+2},y_{2r}(t_{n+2}))\\
h{f}^{\pm}_{2r}(t_{n+1},y_{2r}(t_{n+1}))\\
h{f}^{\pm}_{2r}(t_{n},y_{2r}(t_{n}))
\end{array}\right),
\end{equation*}
and
\begin{equation*}
  \left(
\begin{tabular}{l|lllll}
0&1&$\frac{55}{24}$&$\frac{-59}{24}$&$\frac{37}{24}$&$\frac{-9}{24}$\\
 \cline{1-6}
0&1&$\frac{55}{24}$&$\frac{-59}{24}$&$\frac{37}{24}$&$\frac{-9}{24}$\\
1&0&0&0&0&0\\
0&0&1&0&0&0\\
0&0&0&1&0&0\\
0&0&0&0&1&0
\end{tabular}\right).
\end{equation*}

similarly, for $k=5$ we obtain
\begin{equation*}
{y}^{\mp[n-1]}_{1r}=\left(\begin{array}{c}
{y}^{\mp}_{1r}(t_{n+4})\\
h{f}^{\mp}_{1r}(t_{n+4},y_{1r}(t_{n+4}))\\
h{f}^{\mp}_{1r}(t_{n+3},y_{1r}(t_{n+3}))\\
h{f}^{\mp}_{1r}(t_{n+2},y_{1r}(t_{n+2}))\\
h{f}^{\mp}_{1r}(t_{n+1},y_{1r}(t_{n+1}))\\
h{f}^{\mp}_{1r}(t_{n},y_{1r}(t_{n}))
\end{array}\right),\quad
{y}^{\mp[n-1]}_{2r}=\left(\begin{array}{c}
{y}^{\mp}_{2r}(t_{n+4})\\
h{f}^{\pm}_{2r}(t_{n+4},y_{2r}(t_{n+4}))\\
h{f}^{\pm}_{2r}(t_{n+3},y_{2r}(t_{n+3}))\\
h{f}^{\pm}_{2r}(t_{n+2},y_{2r}(t_{n+2}))\\
h{f}^{\pm}_{2r}(t_{n+1},y_{2r}(t_{n+1}))\\
h{f}^{\pm}_{2r}(t_{n},y_{2r}(t_{n}))
\end{array}\right),
\end{equation*}
and

\begin{equation*}
\left(
\begin{tabular}{l|llllll}
0&1&$\frac{1901}{720}$&$\frac{-2774}{720}$&$\frac{2616}{720}$&$\frac{-1274}{720}$&$\frac{251}{720}$\\
 \cline{1-7}
0&1&$\frac{1901}{720}$&$\frac{-2774}{720}$&$\frac{2616}{720}$&$\frac{-1274}{720}$&$\frac{251}{720}$\\
1&0&0&0&0&0&0\\
0&0&1&0&0&0&0\\
0&0&0&1&0&0&0\\
0&0&0&0&1&0&0\\
0&0&0&0&0&1&0
\end{tabular}\right).
\end{equation*}

\section{Convergence, consistency and stability}
To address the convergence of the presented FGLMs we consider the
numerical solutions
${y}_{1}(t_{n+j};r)=[{y}^-_{1}(t_{n+j};r),{y}^+_{1}(t_{n+j};r)]$ and
${y}_{2}(t_{n+j};r)=[{y}^-_{2}(t_{n+j};r),{y}^+_{2}(t_{n+j};r)]$ and
the corresponding exact solutions
$\mathbf{Y}_{1}(t_{n+j};r)=[\mathbf{Y}^-_{1}(t_{n+j};r),\mathbf{Y}^+_{1}(t_{n+j};r)]$
and
$\mathbf{Y}_{2}(t_{n+j};r)=[\mathbf{Y}^-_{2}(t_{n+j};r),\mathbf{Y}^+_{2}(t_{n+j};r)]$
under (i) and (ii)-differentiability, respectively. The local
truncation errors (LTEs) of the FGLMs under strongly generalized
differentiability are defined by
\begin{equation}\label{eq3.23}
\begin{array}{c}
{\Psi}_{1}(t_{n+k};r)=\sum_{j=0}^{k}r_{j}y_{1}(t_{n+j};r)-h\psi_{f_{1}}\big(y_1(t_{n+k};r),\cdots,y_1(t_n;r)\big),\\
{\Psi}_{2}(t_{n+k};r)=\sum_{j=0}^{k}r_{j}y_{2}(t_{n+j};r)-h\psi_{f_{2}}\big(y_2(t_{n+k};r),\cdots,y_2(t_n;r)\big),\\
\end{array}
\end{equation}
where $r_{k}=-r_{k-1}=1$ and $r_{j}=0$ for $j=0,1,\ldots,k-2$ and
\begin{eqnarray*}
\psi_{f_{1}}\big(y_1(t_{n+k};r),\cdots,y_1(t_n;r)&=&\sum_{j=0}^{k-1}\beta_jf_1(t_{n+j},y_1(t_{n+j};r))\\
\psi_{f_{2}}\big(y_2(t_{n+k};r),\cdots,y_2(t_n;r)&=&\sum_{j=0}^{k-1}\beta_jf_2(t_{n+j},y_2(t_{n+j};r))
\end{eqnarray*}
Consistency and stability are two essential conditions for
convergent.
\begin{defn}
A Fuzzy GLM form of k-step method under generalized
differentiability is said to be consistent if for all fuzzy initial
value problems, the residual ${\Psi}_{1}(t_{n+k};r)$ and
${\Psi}_{2}(t_{n+k};r)$ defined by (\ref{eq3.23}) satisfies
\begin{eqnarray*}
\lim_{h\rightarrow0}\frac{1}{h}{\Psi}_{1}(t_{n+k};r)=0,\\
\lim_{h\rightarrow0}\frac{1}{h}{\Psi}_{2}(t_{n+k};r)=0.
\end{eqnarray*}
\end{defn}
\begin{defn}
A Fuzzy GLM is stable if the minimal polynomial of coefficient matrix $V$ has no
zeros greater than 1 and all zeros equal to 1 are simple, in other words it satisfies the root condition.
\end{defn}
To verify the stability of the given Fuzzy GLMs under generalized
differentiability in section \ref{sec3} we found the minimal
polynomial $p_{k}(w)$ of the coefficient matrix $V$ for $k=4,5$:
\begin{eqnarray*}
p_{k}(w)=w^{k}(w-1),\quad k=4,5,
\end{eqnarray*}
which simply satisfies the root condition and the corresponding Fuzzy GLMs are stable.
\section{Numerical results}\label{sec5}
In this section, we report among many test problems an example to
show the numerical results of FGLMs for solving fuzzy differential
equations under strongly generalized differentiability. We utilize
the FGLMs ($k=4,5$) presented in section \ref{sec3}. The absolute
error numerical results concerning the order of convergence is
provided. We can estimate the order of convergence $p$ by evaluation
of the fraction $\frac{E(h/2)}{E(h)}= O(\frac{1}{2^p})$.

\begin{test}\label{test6.1}
(Bede \cite{Bede1}) Consider the following fuzzy initial value
problem
\begin{equation}\label{FIVP1}
  y'=-y+e^{-t}(-1,0,1),\qquad y_0=(-1,0,1).
\end{equation}
The system of ODEs corresponding to (i)-differentiability is given by
\begin{equation*}
  \left\{\begin{array}{l}
           (y^-)'= -y^+-e^{-t},\\
           (y^1)'= -y^1,\\
           (y^+)' = -y^-+e^{-t}, \\
           y_0=(-1,0,1).
         \end{array}\right.
\end{equation*}
The analytical solution under (i)-differentiability is
\begin{eqnarray}
  Y_1^-(t;r) &=& (1-r)(\frac{1}{2}e^{-t}-\frac{3}{2}e^t) \nonumber\\
  Y_1^+(t;r) &=& (1-r)(\frac{3}{2}e^t-\frac{1}{2}e^{-t}) \nonumber
\end{eqnarray}
Similarly, the system of ODEs corresponding to (ii)-differentiability is given by
\begin{equation*}
  \left\{\begin{array}{l}
           (y^-)'= -y^-+e^{-t},\\
           (y^1)'= -y^1,\\
           (y^+)' = -y^+-e^{-t}, \\
           y_0=(-1,0,1),
         \end{array}\right.
\end{equation*}
and the analytical solution under (ii)-differentiability is
\begin{eqnarray}
  Y_2^-(t;r) &=& (-1+r)(1-t)\exp(-t) \nonumber\\
  Y_2^+(t;r) &=& (1-r)(1-t)\exp(-t). \nonumber
\end{eqnarray}
We demonstrate the numerical solution of FIVP \eqref{FIVP1} in the
interval $[0,2]$. The (i) and (ii)-exact and approximate solutions,
resulted by FGLMs for $k=4$ and $k=5$, are presented in Tables
\ref{Tab6.1.1} and \ref{Tab6.1.2} at $t=2$ with $N=20$ and
$h=\frac{T-t_0}{N}$. Moreover, the results for their convergence
provided in Tables \ref{Tab6.1.3} and \ref{Tab6.1.4}.
\end{test}
\begin{table}[!htp]
  \centering
 \tiny
\begin{tabular}{cccc}
  \hline
  $r$ & $y_{1r}$ & $Y_{1r}$ & $E_{1r}$ \\[0.5mm] \hline\\[-1.5mm]
 0      &   [-1.101531E1, 1.101531E1]   &   [-1.101592E1, 1.101592E1] &     6.024101E-4\\[0.5mm]
0.1     &   [-9.913783E0, 9.913783E0]   &   [-9.914325E0, 9.914325E0] &     5.421691E-4\\[0.5mm]
0.2     &   [-8.812251E0, 8.812251E0]   &   [-8.812733E0, 8.812733E0] &     4.819281E-4\\[0.5mm]
0.3     &   [-7.710720E0, 7.710720E0]   &   [-7.711142E0, 7.711142E0] &     4.216871E-4\\[0.5mm]
0.4     &   [-6.609188E0, 6.609188E0]   &   [-6.609550E0, 6.609550E0] &     3.614461E-4\\[0.5mm]
0.5     &   [-5.507657E0, 5.507657E0]   &   [-5.507958E0, 5.507958E0] &     3.012050E-4\\[0.5mm]
0.6     &   [-4.406126E0, 4.406126E0]   &   [-4.406367E0, 4.406367E0] &     2.409640E-4\\[0.5mm]
0.7     &   [-3.304594E0, 3.304594E0]   &   [-3.304775E0, 3.304775E0] &     1.807230E-4\\[0.5mm]
0.8     &   [-2.203063E0, 2.203063E0]   &   [-2.203183E0, 2.203183E0] &     1.204820E-4\\[0.5mm]
0.9     &   [-1.101531E0, 1.101531E0]   &   [-1.101592E0, 1.101592E0] &     6.024101E-5\\[0.5mm]
1.0     &             [0,0]             &             [0,0]          &  0   \\[0.5mm] \hline
 \\
  (a)
\end{tabular}
\begin{tabular}{cccc}
  \hline
  $r$ & $y_{2r}$ & $Y_{2r}$ & $E_{2r}$\\[0.5mm] \hline\\[-1.5mm]
 0      &   [1.352883E-1, -1.352883E-1] &   [1.353353E-1, -1.353353E-1] &   4.699417E-5\\[0.5mm]
0.1     &   [1.217595E-1, -1.217595E-1] &   [1.218018E-1, -1.218018E-1] &   4.229476E-5\\[0.5mm]
0.2     &   [1.082306E-1, -1.082306E-1] &   [1.082682E-1, -1.082682E-1] &   3.759534E-5\\[0.5mm]
0.3     &   [9.470180E-2, -9.470180E-2] &   [9.473470E-2, -9.473470E-2] &   3.289592E-5\\[0.5mm]
0.4     &   [8.117297E-2, -8.117297E-2] &   [8.120117E-2, -8.120117E-2] &   2.819650E-5\\[0.5mm]
0.5     &   [6.764414E-2, -6.764414E-2] &   [6.766764E-2, -6.766764E-2] &   2.349709E-5\\[0.5mm]
0.6     &   [5.411532E-2, -5.411532E-2] &   [5.413411E-2, -5.413411E-2] &   1.879767E-5\\[0.5mm]
0.7     &   [4.058649E-2, -4.058649E-2] &   [4.060058E-2, -4.060058E-2] &   1.409825E-5\\[0.5mm]
0.8     &   [2.705766E-2, -2.705766E-2] &   [2.706706E-2, -2.706706E-2] &   9.398835E-6\\[0.5mm]
0.9     &   [1.352883E-2, -1.352883E-2] &   [1.353353E-2, -1.353353E-2] &   4.699417E-6\\[0.5mm]
   1.0 &              [0,0]             &             [0,0]             &           0   \\[0.5mm] \hline
   \\
(b)
\end{tabular}
  \caption{\scriptsize(a)Approximate solution of the FGLM ($k=4$) $y_{1r}=[y_{1r}^-,y_{1r}^+]$, exact solution $Y_{1r}=[Y_{1r}^-,Y_{1r}^+]$ and absolute error $E_{1r}$ under (i)-differentiability, (b)Approximate solution of the FGLM ($k=4$) $y_{2r}=[y_{2r}^-,y_{2r}^+]$, exact solution $Y_{2r}=[Y_{2r}^-,Y_{2r}^+]$ and absolute error $E_{2r}$ under (ii)-differentiability, Test \ref{test6.1}.}\label{Tab6.1.1}
\end{table}

\begin{table}[!htp]
  \centering
 \tiny
\begin{tabular}{ccccccc}
  \hline
  $r$ & $y_{1r}$ & $Y_{1r}$ & $E_{1r}$ \\[0.5mm] \hline\\[-1.5mm]
 0      &   [-1.101587E1, 1.101587E1]   &   [-1.101592E1, 1.101592E1]   &   4.451187E-5\\[0.5mm]
0.1     &   [-9.914285E0, 9.914285E0]   &   [-9.914325E0, 9.914325E0]   &   4.006069E-5\\[0.5mm]
0.2     &   [-8.812698E0, 8.812698E0]   &   [-8.812733E0, 8.812733E0]   &   3.560950E-5\\[0.5mm]
0.3     &   [-7.711110E0, 7.711110E0]   &   [-7.711142E0, 7.711142E0]   &   3.115831E-5\\[0.5mm]
0.4     &   [-6.609523E0, 6.609523E0]    &  [-6.609550E0, 6.609550E0]   &   2.670712E-5\\[0.5mm]
0.5     &   [-5.507936E0, 5.507936E0]   &   [-5.507958E0, 5.507958E0]   &   2.225594E-5\\[0.5mm]
0.6     &   [-4.406349E0, 4.406349E0]   &   [-4.406367E0, 4.406367E0]   &   1.780475E-5\\[0.5mm]
0.7     &   [-3.304762E0, 3.304762E0]   &   [-3.304775E0, 3.304775E0]   &   1.335356E-5\\[0.5mm]
0.8     &   [-2.203174E0, 2.203174E0]   &   [-2.203183E0, 2.203183E0]   &   8.902375E-6\\[0.5mm]
0.9     &   [-1.101587E0, 1.101587E0]   &   [-1.101592E0, 1.101592E0]   &   4.451187E-6\\[0.5mm]
1.0     &             [0,0]             &             [0,0]             &      0       \\[0.5mm] \hline
 \\
  (a)
\end{tabular}
\begin{tabular}{cccc}
  \hline
  $r$ & $y_{2r}$ & $Y_{2r}$ & $E_{2r}$ \\[0.5mm] \hline\\[-1.5mm]
 0      &   [1.353406E-1, -1.353406E-1] &   [1.353353E-1, -1.353353E-1] &       5.270043E-6 \\[0.5mm]
0.1     &   [1.218065E-1, -1.218065E-1] &   [1.218018E-1, -1.218018E-1] &       4.743039E-6 \\[0.5mm]
0.2     &   [1.082724E-1, -1.082724E-1] &   [1.082682E-1, -1.082682E-1] &       4.216035E-6 \\[0.5mm]
0.3     &   [9.473839E-2, -9.473839E-2] &   [9.473470E-2, -9.473470E-2] &       3.689030E-6 \\[0.5mm]
0.4     &   [8.120433E-2, -8.120433E-2] &   [8.120117E-2, -8.120117E-2] &       3.162026E-6 \\[0.5mm]
0.5     &   [6.767028E-2, -6.767028E-2] &   [6.766764E-2, -6.766764E-2] &       2.635022E-6 \\[0.5mm]
0.6     &   [5.413622E-2, -5.413622E-2] &   [5.413411E-2, -5.413411E-2] &       2.108017E-6 \\[0.5mm]
0.7     &   [4.060217E-2, -4.060217E-2] &   [4.060058E-2, -4.060058E-2] &       1.581013E-6 \\[0.5mm]
0.8     &   [2.706811E-2, -2.706811E-2] &   [2.706706E-2, -2.706706E-2] &       1.054009E-6 \\[0.5mm]
0.9     &   [1.353406E-2, -1.353406E-2] &   [1.353353E-2, -1.353353E-2] &       5.270043E-7 \\[0.5mm]
1.0     &             [0,0]             &             [0,0]            &             0      \\[0.5mm] \hline
   \\
(b)
\end{tabular}
  \caption{\scriptsize(a)Approximate solution of the FGLM ($k=5$), $y_{1r}=[y_{1r}^-,y_{1r}^+]$, exact solution $Y_{1r}=[Y_{1r}^-,Y_{1r}^+]$ and absolute error $E_{1r}$  under (i)-differentiability, (b)Approximate solution of the FGLM ($k=5$), $y_{2r}=[y_{2r}^-,y_{2r}^+]$, exact solution $Y_{2r}=[Y_{2r}^-,Y_{2r}^+]$ and absolute error $E_{2r}$ under (ii)-differentiability, Test \ref{test6.1}.}\label{Tab6.1.2}
\end{table}

\begin{table}[!htp]
  \centering
  \tiny
  \begin{tabular}{cccccccccccc}
    \hline
    $r$&$h_{i}$&\quad&\quad&\quad&\quad& $E_{2r}(h_{i})$&\quad&\quad&\quad&\quad& $p$\\[1mm] \hline
     0.2& $\frac{1}{10}$  & \quad&\quad&\quad&\quad& 3.759533862475462E-5& \quad&\quad&\quad&\quad&                 \quad\\[1mm]
   \quad& $\frac{1}{20}$  & \quad&\quad&\quad&\quad& 2.360816361970941E-6& \quad&\quad&\quad&\quad&  3.993196066118324E0\\[1mm]
   \quad& $\frac{1}{40}$  & \quad&\quad&\quad&\quad& 1.475860954835984E-7& \quad&\quad&\quad&\quad&  3.999657112312050E0\\[1mm]
   \quad& $\frac{1}{80}$  & \quad&\quad&\quad&\quad& 9.220812974275461E-9& \quad&\quad&\quad&\quad&  4.000519042265663E0\\[1mm]
    \hline
    0.4  &$\frac{1}{10}$  & \quad&\quad&\quad&\quad& 2.819650396852780E-5& \quad&\quad&\quad&\quad&                 \quad\\[1mm]
    \quad&$\frac{1}{20}$  & \quad&\quad&\quad&\quad& 1.770612271481675E-6& \quad&\quad&\quad&\quad&  3.993196066113545E0\\[1mm]
    \quad&$\frac{1}{40}$  & \quad&\quad&\quad&\quad& 1.106895716057599E-7& \quad&\quad&\quad&\quad&  3.999657112405316E0\\[1mm]
    \quad&$\frac{1}{80}$  & \quad&\quad&\quad&\quad& 6.915609765401065E-9& \quad&\quad&\quad&\quad&  4.000519034937461E0\\[1mm]
    \hline
    0.6  &$\frac{1}{10}$  & \quad&\quad&\quad&\quad& 1.879766931239119E-5& \quad&\quad&\quad&\quad&                 \quad\\[1mm]
    \quad&$\frac{1}{20}$  & \quad&\quad&\quad&\quad& 1.180408180992409E-6& \quad&\quad&\quad&\quad&  3.993196066110909E0\\[1mm]
    \quad&$\frac{1}{40}$  & \quad&\quad&\quad&\quad& 7.379304775567697E-8& \quad&\quad&\quad&\quad&  3.999657112049211E0\\[1mm]
    \quad&$\frac{1}{80}$  & \quad&\quad&\quad&\quad& 4.610406514893306E-9& \quad&\quad&\quad&\quad&  4.000519033851667E0\\[1mm]
    \hline
    0.8  &$\frac{1}{10}$  & \quad&\quad&\quad&\quad& 9.398834656195593E-6& \quad&\quad&\quad&\quad&                 \quad\\[1mm]
    \quad&$\frac{1}{20}$  & \quad&\quad&\quad&\quad& 5.902040904962047E-7& \quad&\quad&\quad&\quad&  3.993196066110909E0\\[1mm]
    \quad&$\frac{1}{40}$  & \quad&\quad&\quad&\quad& 3.689652387783848E-8& \quad&\quad&\quad&\quad&  3.999657112049211E0\\[1mm]
    \quad&$\frac{1}{80}$  & \quad&\quad&\quad&\quad& 2.305203257446653E-9& \quad&\quad&\quad&\quad&  4.000519033851667E0\\[1mm]
    \hline
  \end{tabular}
  \caption{\scriptsize Convergence of the FGLM ($k=4$)under (ii)-differentiability}\label{Tab6.1.3}
\end{table}
\begin{table}[!htp]
  \centering
  \tiny
  \begin{tabular}{cccccccccccc}
    \hline
    $r$&$h_{i}$&\quad&\quad&\quad&\quad& $E_{2r}(h_{i})$&\quad&\quad&\quad&\quad& $p$\\ [.5mm]\hline\\[-1.5mm]
 0.2 & $\frac{1}{10}$  & \quad&\quad&\quad&\quad& 4.216034534335056E-6   & \quad&\quad&\quad&\quad&                 \quad\\[1mm]
\quad& $\frac{1}{20}$  & \quad&\quad&\quad&\quad& 1.336319765954386E-7   & \quad&\quad&\quad&\quad&  4.979549509841708E0\\[1mm]
\quad& $\frac{1}{40}$  & \quad&\quad&\quad&\quad& 4.185707641601866E-9   & \quad&\quad&\quad&\quad&  4.996649911789974E0\\[1mm]
\quad& $\frac{1}{80}$  & \quad&\quad&\quad&\quad& 1.308334413030465E-10  & \quad&\quad&\quad&\quad&  4.999668298467584E0\\[1mm]
\quad& $\frac{1}{160}$ & \quad&\quad&\quad&\quad& 4.088021587911328E-12  & \quad&\quad&\quad&\quad&  5.000184718796247E0\\[1mm]
    \hline
 0.4 & $\frac{1}{10}$  & \quad&\quad&\quad&\quad& 3.162025900754761E-6  & \quad&\quad&\quad&\quad&                 \quad\\[1mm]
\quad& $\frac{1}{20}$  & \quad&\quad&\quad&\quad& 1.002239824188234E-7  & \quad&\quad&\quad&\quad&  4.979549510242824E0\\[1mm]
\quad& $\frac{1}{40}$  & \quad&\quad&\quad&\quad& 3.139280738140293E-9  & \quad&\quad&\quad&\quad&  4.996649908201587E0\\[1mm]
\quad& $\frac{1}{80}$  & \quad&\quad&\quad&\quad& 9.812499424111110E-11 & \quad&\quad&\quad&\quad&  4.999669576905354E0\\[1mm]
\quad& $\frac{1}{160}$ & \quad&\quad&\quad&\quad& 3.065672715685253E-12 & \quad&\quad&\quad&\quad&  5.000345072764134E0\\[1mm]
\hline
 0.6 & $\frac{1}{10}$  & \quad&\quad&\quad&\quad& 2.108017267167528E-6  & \quad&\quad&\quad&\quad&                  \quad\\[.5mm]
\quad& $\frac{1}{20}$  & \quad&\quad&\quad&\quad& 6.681598831159707E-8  & \quad&\quad&\quad&\quad&  4.979549509542057E0\\[.5mm]
\quad& $\frac{1}{40}$  & \quad&\quad&\quad&\quad& 2.092853827739827E-9  & \quad&\quad&\quad&\quad&  4.996649907306344E0\\[.5mm]
\quad& $\frac{1}{80}$  & \quad&\quad&\quad&\quad& 6.541663044590251E-11 & \quad&\quad&\quad&\quad&  4.999670292639665E0\\[.5mm]
\quad& $\frac{1}{160}$ & \quad&\quad&\quad&\quad& 2.043996916167856E-12 & \quad&\quad&\quad&\quad&  5.000192524602108E0\\[.5mm]
\hline
 0.8 & $\frac{1}{10}$  & \quad&\quad&\quad&\quad& 1.054008633583764E-6  & \quad&\quad&\quad&\quad&                 \quad\\[1mm]
\quad& $\frac{1}{20}$  & \quad&\quad&\quad&\quad& 3.340799415579854E-8  & \quad&\quad&\quad&\quad&  4.979549509542057E0\\[1mm]
\quad& $\frac{1}{40}$  & \quad&\quad&\quad&\quad& 1.046426913869913E-9  & \quad&\quad&\quad&\quad&  4.996649907306344E0\\[1mm]
\quad& $\frac{1}{80}$  & \quad&\quad&\quad&\quad& 3.270831522295126E-11 & \quad&\quad&\quad&\quad&  4.999670292639665E0\\[1mm]
\quad& $\frac{1}{80}$  & \quad&\quad&\quad&\quad& 1.021998458083928E-12 & \quad&\quad&\quad&\quad&  5.000192524602108E0\\[1mm]
\hline
  \end{tabular}
  \caption{\scriptsize Convergence of the FGLM ($k=5$)under (ii)-differentiability}\label{Tab6.1.4}
\end{table}

From Tables \ref{Tab6.1.3} and \ref{Tab6.1.4}, it follows that the
Fuzzy GLMs of $4$-step methods under strongly generalized
differentiability have convergence order 4 and the Fuzzy GLMs form
of 5-step methods have convergence order 5.

\section{Conclusion}
In this paper we have developed the linear multistep methods
(Adams-Bashforth methods) in the framework of general linear methods
for solving fuzzy differential equations under strongly generalized
differentiability. We have shown the consistency, stability, and
convergence of the new FGLM formulation. The general framework of
FGLMs will be studied in the forthcoming paper.

\end{document}